\documentclass[12pt]{amsart}
\usepackage{amssymb}
\usepackage{amsmath}
\usepackage{longtable}
\newcommand\g{{\mathfrak g}}
\renewcommand\b{{\mathfrak b}}
\newcommand\h{{\mathfrak h}}

\newcommand{\p}{\mathfrak{p}}
\newcommand\gl{\mathfrak{gl}}
\renewcommand{\a}{\mathfrak{a}}

\renewcommand\l{\mathfrak l}

\newcommand\z{\mathfrak z}
\newcommand\s{\mathfrak s}
\renewcommand{\t}{\mathfrak{t}}
\newcommand\codim{\operatorname{codim}}

\newcommand\C{\mathbb C}
\newcommand\Q{\mathbb Q}
\newcommand\A{\mathfrak A}
\newcommand\X{\mathfrak X}
\newcommand\R{\mathbb R}

\newcommand\N{\mathbb N}
\newcommand\Z{\mathbb Z}

\newcommand\SL{\operatorname{SL}}
\newcommand\Sp{\operatorname{Sp}}

\newcommand{\Ad}{\mathop{\rm Ad}\nolimits}

\newcommand{\rank}{\mathop{\rm rk}\nolimits}

\newcommand\Int{\mathop{\rm Int}\nolimits}
\newcommand\Aut{\mathop{\rm Aut}\nolimits}

\newtheorem{Thm}{Theorem}[section]
\newtheorem{Prop}[Thm]{Proposition}
\newtheorem{Cor}[Thm]{Corollary}
\newtheorem{Lem}[Thm]{Lemma}
\theoremstyle{definition}
\newtheorem{Ex}[Thm]{Example}
\newtheorem{defi}[Thm]{Definition}

\newtheorem{Prob}[Thm]{Problem}
\unitlength=1mm \numberwithin{equation}{section} \oddsidemargin=0cm
\author{Ivan V. Losev}
\title{Embeddings of homogeneous spaces into irreducible modules}
\begin{document}
\begin{abstract}
Let $G$ be a connected reductive group. We find a necessary and
sufficient condition for a quasiaffine homogeneous space of $G/H$ to
be embeddable into an irreducible $G$-module.  If $H$ is reductive
we also find a necessary and sufficient condition for a closed
embedding of $G/H$ into an irreducible module to exist. These
conditions are stated in terms of the group of {\it central
automorphisms} of $G/H$.
\end{abstract}
\maketitle
\section{Introduction}
The base field is the field $\C$ of complex numbers.
Throughout the paper $G$ denotes a connected reductive
algebraic group, $B$ its Borel subgroup and $T$ a maximal
torus of $B$.

The celebrated theorem of Chevalley states that any homogeneous
space can be embedded (as a locally-closed subvariety) into the
projectivization of a $G$-module. If $H$ is an observable subgroup
of $G$, that is, the homogeneous space $G/H$ is quasiaffine, then
$G/H$ can be embedded even into a $G$-module itself, see, for
example, \cite{VP}, Theorem 1.6. So it is natural to pose the
following
\begin{Prob}\label{Prob:1.0}
Describe all observable subgroups $H$ such that $G/H$ can be embedded
into an {\it irreducible} $G$-module.
\end{Prob}

To state the answer to that problem we need the definition of a
central automorphism of a $G$-variety. Let $X$ be an irreducible
$G$-variety.  A subspace $\C(X)^{(B)}_{\lambda}\subset\C(X)$
consisting of all $B$-semiinvariant functions of weight $\lambda\in
\X(B)$
 is stable
under every $G$-equivariant automorphism of $X$. The following
definition belongs to Knop, \cite{Knop8}.

\begin{defi}\label{Def:1.1}
A $G$-equivariant automorphism of  $X$ is called {\it central} if it
acts on  $\C(X)^{(B)}_{\lambda}$ by the multiplication by a constant
for any  weight $\lambda$.
\end{defi}

We denote the group of central automorphisms of $X$ by $\A_G(X)$. We
write $\A_{G,H}$ instead of $\A_{G}(G/H)$. It was shown by Knop,
\cite{Knop8}, Section 5, that $\A_{G,H}$ is an algebraic
quasi-torus, that is a closed subgroup of an algebraic torus.

\begin{Thm}\label{Thm:1.2}
Let $H$ be an observable subgroup of $G$. Then the following conditions are
equivalent:
\begin{itemize}
\item[(a)] $G/H$ can be embedded into an irreducible $G$-module.
\item[(b)] $\A_{G,H}$ is a finite cyclic group or a one-dimensional
torus.
\end{itemize}
\end{Thm}

For a given subgroup $H\subset G$ the group $\A_{G,H}$ can be
computed using techniques from \cite{weight}. Namely, $\A_{G,H}$ is
the quotient of the  {\it weight lattice} of $G/H$ by the {\it root
lattice} of $G/H$. An algorithm for computing the weight lattice is
the main result of \cite{weight}. The computation of the root
lattice can be reduced to that of the weight lattice by using
\cite{weight}, Proposition 5.2.1.

If $H$ is a reductive subgroup of $G$ or, equivalently, $G/H$ is
affine, then one may ask whether there exists a {\it closed}
embedding of $G/H$ into an irreducible $G$-module. Here is an
answer.

\begin{Thm}\label{Thm:1.3}
Let $H$ be a reductive subgroup of $G$. Then the following
conditions are equivalent.
\begin{itemize}
\item[(a)] There is a closed equivariant embedding of $G/H$ into a
irreducible $G$-module.
\item[(b)] $\A_{G,H}$ is a finite cyclic group.
\end{itemize}
\end{Thm}

We prove Theorems \ref{Thm:1.2},\ref{Thm:1.3} in Sections
\ref{SECTION_Thm1},\ref{SECTION_Thm2}. In Section
\ref{SECTION_Examples} we present some examples of applications of
our theorems.

{\bf Acknowledgements.} I thank I.V. Arzhantsev, who communicated
this problem to me, found a gap in the proof in the previous version
of the paper and made some other useful remarks. The paper was
partially written during my stay  in the Fourier University,
Grenoble, in June, 2006. I express my gratitude to this institution
and especially to Prof. M. Brion for hospitality.

\section{Notation and conventions}
\setlongtables
\begin{longtable}{p{4.5cm} p{11cm}}
$A^{(B)}_\mu$& the subspace of all $B$-semiinvariant functions
of weight $\mu$ in a $G$-algebra $A$, where $G$ is a connected
reductive group.
\\ $[\g,\g]$& the commutant of a Lie algebra $\g$.
\\ $G^{\circ}$& the connected component of unit of an algebraic
group $G$.
\\ $R_u(G)$& the unipotent radical of an algebraic group $G$.
\\ $G_x$& the stabilizer of a point $x\in X$ under an action
$G:X$.
\\
$\Int(\g)$& the group of  inner
automorphisms of a Lie algebra $\g$.
\\ $N_G(H)$& the normalizer of a subgroup $H$ in a group $G$
\\ $V^\g$& $=\{v\in V| \g v=0\}$, where $\g$ is a Lie algebra and
$V$ is a $\g$-module.
\\ $V(\mu)$& the irreducible module with highest weight $\mu$
over a reductive algebraic group  or a reductive Lie algebra.\\
$\X(G)$& the character lattice of an algebraic group $G$.\\   $X^G$& the
fixed-point set
for an action of $G$ on $X$.\\
$\#X$& the cardinality of a  set $X$.
\\$Z(G)$ (resp., $\z(\g)$)& the center of an algebraic group $G$ (resp., of a Lie algebra $\g$).
\\   $Z_G(\h)$ (resp., $\z_\g(\h)$)& the centralizer  of a
subalgebra $\h\subset \g$  in an algebraic group $G$ (resp., in its Lie algebra $\g$).
\\ $\lambda^*$& the  highest weight  dual to $\lambda$.
\end{longtable}

If an algebraic
group is denoted by a capital Latin letter, then we denote its Lie algebra
by the corresponding small frakture letter, for example, $\widehat{\h}$
denotes the Lie algebra of $\widehat{H}$. All topological terms refer to the
Zariski topology.

\section{Proof of Theorem~\ref{Thm:1.2}}\label{SECTION_Thm1}
First, we fix some notation and recall some definitions from the
theory of algebraic transformation groups.

In this section $H$ denotes an observable subgroup of $G$. The group
of $G$-equivariant automorphisms of $G/H$ is identified with
$N_G(H)/H$. We consider $\A_{G,H}$ as a subgroup in $N_G(H)/H$.
Denote by $H^{sat}$ the inverse image of $\A_{G,H}$ in $N_G(H)$.

Let $X$ be an irreducible $G$-variety. An element $\lambda\in \X(T)$
is said to be a {\it weight of $X$} if $\C(X)_\lambda^{(B)}\neq 0$.
Clearly, all weights of $X$ form a subgroup of $\X(T)$ called the {\it weight
lattice} of $X$ and  denoted by $\X_{G,X}$. The rank of the weight lattice is called the {\it rank} of
$X$ and is  denoted by $\rank_G(X)$. We put $\a_{G,X}=\X_{G,X}\otimes_{\Z}\C$.
If $X=G/G_0$, then we write $\X_{G,G_0}$ instead of $\X_{G,G/G_0}$. It is easy to
see that the subspace $\a_{G,G/G_0}$ depends only on the pair $(\g,\g_0)$. Thus we write
$\a_{\g,\g_0}$ instead of $\a_{G,G/G_0}$. If $\widehat{G}_0$ is a subgroup of $G$ containing
$G_0$, then there exists a dominant $G$-equivariant morphism $G/G_0\rightarrow G/\widehat{G}_0$ and
thence $\X_{G,\widehat{G}_0}\subset \X_{G,G_0}$.

The codimension of a general $B$-orbit in $X$ is called the {\it
complexity} of $X$ and is denoted by $c_G(X)$. Again, we write
$c_{\g,\g_0}$ instead of $c_G(G/G_0)$. Let us note that
$c_{\g,\widehat{\g}_0}\leqslant c_{\g,\g_0}$ whenever $G_0\subset
\widehat{G}_0$.

Proceed to the proof of Theorem~\ref{Thm:1.2}. The implication (a)$\Rightarrow$(b)
is easy.

\begin{proof}[Proof of (a)$\Rightarrow$(b)]
By the Frobenius reciprocity, there is an $N_G(H)$-equivariant isomorphism
$V(\lambda)^H\cong \C[G/H]^{(B)}_{\lambda^*}$. Clearly, (a) implies that the action
of $N_G(H)/H$ on $V(\lambda)^H$ is effective for some $\lambda$. Now (b) follows
easily from the definition of the subgroup $\A_{G,H}\subset N_G(H)/H$.
\end{proof}

The implication $(b)\Rightarrow (a)$ will follow from the following
\begin{Prop}\label{Prop:2.0}
Suppose $\A_{G,H}$ is a cyclic finite group or a one-dimensional
torus. Then there is a highest weight $\lambda$ such that
$V(\lambda)^H\neq \{0\}$ and the subset
$\bigcap_{\widehat{H}\supsetneq H} V(\lambda)^{\widehat{H}}$ is not
dense in $ V(\lambda)^H$.
\end{Prop}


The scheme of the proof of the proposition is, roughly speaking, as
follows. On the first step  we prove that for an appropriate highest
weight $\lambda$ the complexity $c_{\g,\g_v}$  for a point $v\in
V(\lambda)^H$ in general position coincides with $c_{\g,\h}$. On the
second step we check that one may choose $\lambda$ such that
$\g_v=\h$ for $v\in V(\lambda)^H$ in general position. At last, we
show that $G_v=H$ for  general $v\in V(\lambda)^H$.

We begin with some simple lemmas.


\begin{Lem}\label{Lem:2.1}
$\dim V(\nu)^H\leqslant \dim V(\nu+\mu)^H$ for any highest weights
$\mu,\nu$ such that $V(\mu)^H\neq 0$.
\end{Lem}
\begin{proof}
By the Frobenius reciprocity, $V(\nu)^H\cong \C[G/H]^{(B)}_{\nu^*}$,
$V(\nu+\mu)^H\cong \C[G/H]^{(B)}_{(\nu+\mu)^*}$. The map $f_1\mapsto
ff_1: \C[G/H]^{(B)}_{\nu^*}\hookrightarrow
\C[G/H]^{(B)}_{(\nu+\mu)^*}$ is injective  for any $f\in
\C[G/H]^{(B)}_{\mu^*},f\neq 0$.
\end{proof}

In the sequel we will need some properties of central automorphisms.

\begin{Lem}\label{Lem:2.2}
\begin{enumerate}
\item An element $n\in N_G(H)/H$ is  central  iff it acts trivially
on $K(G/H)^{(B)}$.
\item $\A_{G,H}\subset Z(N_G(H)/H)$.
\end{enumerate}
\end{Lem}
\begin{proof}
Let $X$ be an affine $G$-variety with open $G$-orbit $G/H$. Thanks
to \cite{VP}, Theorem 3.3, to prove assertion 1 it is enough to
check that $n$ acts on $\C[X]^{(B)}_\lambda$ by the multiplication
by a constant for any highest weight $\lambda$ provided $n$ acts
trivially on $\C(G/H)^{B}$. Since $X$ contains a dense $G$-orbit, we
have $\C[X]^G=\C$. It follows from \cite{VP}, Theorem 3.24, that
$\dim \C[X]^{(B)}_\lambda<\infty$. Now our claim is clear.

Assertion 2 follows from \cite{Knop8}, Corollary 5.6.
\end{proof}

The following technical proposition is crucial in the proof of
Proposition \ref{Prop:2.0}.

\begin{Prop}\label{Prop:2.1}
Let $\a_1,\ldots,\a_k$ be proper subspaces of $\a_{\g,\h}$
and $\X_1,\ldots, \X_l$  sublattices of $\X_{G,H}$ such that $p_i:=\# (\X_{G,H}/\X_i), i=\overline{1,l},$
are pairwise different primes.
Put $c:=c_{\g,\h}$.
Then there exists a highest weight $\lambda$ satisfying
condition (1), when $c$ is arbitrary, and conditions (2),(3),
when $c>0$.
\begin{enumerate}
\item $\lambda^*\not\in\bigcup_{i=1}^k \a_i\cup\bigcup_{i=1}^l\X_i$.
\item The codimension of the closure of the subset $Z:=(\bigcup V(\lambda)^{\widehat{\h}})\cap
V(\lambda)^H)\subset V(\lambda)^H$, where the union is taken over
all algebraic subalgebras $\widehat{\h}\subset\g$ such that
$\widehat{\h}\supset\h, c_{\g,\widehat{\h}}< c$, is strictly  bigger
than $2\dim G$.
\item For any $f\in \C(G/H)^B$ there exist $f_1,f_2\in \C[G/H]^{(B)}_{\lambda^*}$
such that $f=\frac{f_1}{f_2}$.
\end{enumerate}
\end{Prop}

\begin{Lem}\label{Lem:3.0}
Let $\a_1,\ldots,\a_k,\X_1,\ldots,\X_l$ be such as in Proposition
\ref{Prop:2.1}. Let $\mu'\in\Psi$ satisfy condition (1). Then there
is $n\in\N$ such that for any $\lambda\in\Psi$ at least one of the
weights $\lambda+\mu',\lambda+2\mu',\ldots,\lambda+n\mu'$ satisfies
condition (1) of Proposition \ref{Prop:2.1}.
\end{Lem}
\begin{proof}
Set $n:=(k+1)p_1\ldots p_l$. The proof is easy.
\end{proof}

\begin{proof}[Proof of Proposition \ref{Prop:2.1}]
Let us choose a norm $|\cdot|$ on the space $\a_{\g,\h}(\R):=\X_{G,H}\otimes_\Z \R$.
It follows from Timashev's theorem, \cite{Timashev}, that the following assertions hold:
\begin{itemize}
\item There exists $A_0\in \R$ such that   $\dim V(\lambda)^{\widehat{\h}}< A_0 |\lambda|^{c-1}$
for any subalgebra $\widehat{\h}\subset \g$ with
$c_{\g,\widehat{\h}}<c$ and  any highest weight $\lambda$.
\item For any $A\in \R$ there exists a highest weight $\lambda$ such that
$\dim V(\lambda)^{H}> A |\lambda|^{c-1}$.
\end{itemize}

Denote by $Y$ the subvariety of
$\coprod_{i=\dim\h}^{\dim\g}\operatorname{Gr}_i(\g)$ consisting of
all subalgebras $\widehat{\h}\subset \g$ containing $\h$.
$Y_0:=\{\widehat{\h} \in Y| c_{\g,\widehat{\h}}<c\}$ is an open
subvariety of $Y$, because $c_{\g,\widehat{\h}}=\min_{g\in G} \dim
\g/(\Ad(g)\b+\widehat{\h})$. Put $V:=V(\lambda)^H$,
$\widetilde{Z}:=\{ (\widehat{\h},v)\in Y_0\times V| v\in
V(\lambda)^{\widehat{\h}}\}$. The latter is a closed subvariety in
$Y_0\times V$
 of dimension at most $\dim Y_0+\max_{\widehat{\h}\in Y_0}
\dim V(\lambda)^{\widehat{\h}}$.

Note that $Z$ is just the image of $\widetilde{Z}$ under the
projection $Y_0\times V\rightarrow V$. Thus if $c>0$, then the
dimension of the closure of $Z$ does not exceed $
A_0|\lambda|^{c-1}+\dim Y_0$.

Note that there exists a highest weight $\lambda_1$ satisfying
condition (3). Indeed, the field $\C(G/H)^B$ is finitely generated
and let $f_1,\ldots, f_s$ be its generators. Analogously to
\cite{VP}, Theorem 3.3, one proves that there are $f_{i1},f_{i2}\in
\C[G/H]^{(B)}_{\nu_i}, i=\overline{1,s},$ such that
$f_i=\frac{f_{i1}}{f_{i2}}$. It is enough to take $\sum_{i=1}^s
\nu_i^*$ for $\lambda_1$. Note that for any highest weight
$\lambda_2$ with $\C[G/H]^{(B)}_{\lambda_2^*}\neq 0$ the highest
weight $\lambda_2+\lambda_1$ also satisfies condition (3).

Note that there is a highest weight $\lambda_2$ satisfying condition
(1) and such that $V(\lambda_2)^H\neq \{0\}$. Indeed, otherwise
$\cup_{i=1}^k \a_i$ contains an open cone in $\a_{\g,\h}$, which is
absurd. So in case $c=0$ we are done. Now suppose $c>0$.

Let $n$ be such as in Lemma \ref{Lem:3.0}. Choose $A>0$ and a
highest weight $\nu$ such that
    $\dim V(\nu)^H> A|\nu|^{c-1}$
and $A|\nu|^{c-1}> A_0(|\nu|+|\lambda_1|+n|\lambda_2|)^{c-1}+\dim
Y_0+2\dim G.$  Further, there is $j\in \{1,\ldots, n\}$ such that
$\lambda:= \nu+\lambda_1+j\lambda_2$ satisfies (1).  For $c>0$ it is
easy to deduce from Lemma~\ref{Lem:2.1} that $\lambda$ satisfies
condition (2). Finally $\lambda$ satisfies condition (3), for it is
of the form $\lambda_1+\lambda_2'$ for some $\lambda_2$ with
$\C[G/H]^{(B)}_{\lambda_2^*}\neq 0$.
\end{proof}

The next proposition is used on the second step  of the proof.

\begin{Prop}\label{Prop:2.2}
The set $\{\a_{\g,\widehat{\h}}|
\widehat{\h}=[\widehat{\h},\widehat{\h}]+R_u(\widehat{\h})+\h,
\widehat{\h}\text{ is algebraic}\}$ is finite.
\end{Prop}
\begin{proof}
Let $\h=\s\oplus R_u(\h), \widehat{\h}=\widehat{\s}\oplus R_u(\widehat{\h})$ be Levi decompositions.
We may assume that $\s\subset \widehat{\s}$. Denote by $\widehat{H},\widehat{S}$
the connected subgroups of $G$ corresponding to $\widehat{\h},\widehat{\s}$.
By the Weisfeller theorem, see \cite{Weisfeller}, there is a parabpolic subgroup
$P\subset G$ and a Levi subgroup $L\subset P$ such that $\widehat{S}\subset L$,
$R_u(\widehat{H})\subset R_u(P)$.
Conjugating $\h,\widehat{\h}$ by an element of $G$, we may assume that $T\subset L$ and that
$P$ is opposite to $B$. By Panyushev's theorem, \cite{Panyushev2},
$$\a_{\g,\widehat{\h}}=\a_{L,L*_{\widehat{S}}(R_u(\p)/R_u(\widehat{\h}))}.$$
There is an inclusion of $\widehat{S}$-modules
$R_u(\p)/R_u(\widehat{\h})\hookrightarrow \g/\widehat{\s}$. So the
set $\{(L,R_u(\p)/R_u(\widehat{\h}))\}$ is finite. It remains to
check that $\widehat{\s}$ belongs only to finitely many
$\Int(\l)$-conjugacy classes. The following well-known lemma (which
stems, for example, from \cite{Vinberg2}, Proposition 3) allows us
to replace $\Int(\l)$-conjugacy in the previous statement by
$\Int(\g)$-conjugacy.

\begin{Lem}\label{Lem:2.3}
Let $\g_0$ be a reductive subalgebra of $\g$ and $\g_1$ a reductive
subalgebra of $\g_0$. The set of subalgebras of $\g_0$, that are $\Int(\g)$-conjugate
to $\g_1$, decomposes into finitely many classes of $\Int(\g_0)$-conjugacy.
\end{Lem}

The equality $\widehat{\h}=[\widehat{\h},\widehat{\h}]+R_u(\widehat{\h})+\h$
is equivalent to $\widehat{\s}=[\widehat{\s},\widehat{\s}]+\s$. Therefore the statement
on the finiteness of the set of $\Int(\g)$-conjugacy classes stems from the following
lemma
\begin{Lem}\label{Lem:2.4}
Let $\s$ be a reductive subalgebra of $\g$. The set of
$\Int(\g)$-conjugacy of reductive subalgebras
$\widehat{\s}\subset\g$ such that
$\widehat{\s}=[\widehat{\s},\widehat{\s}]+\s$ is finite.
\end{Lem}
\begin{proof}[Proof of Lemma~\ref{Lem:2.4}]
We may replace $\s$ with its Cartan subalgebra and assume that
$\s\subset \t$. In this case the proof is in three steps.

{\it Step 1.} Here we show that the set of subspaces of $\t$, that
are Cartan subalgebras of semisimple subalgebras of $\g$, is finite.
Note that there are finitely many conjugacy classes of semisimple
subalgebras of $\g$. Indeed, for $\g=\gl_n$ this is a consequence of
the highest weight theory and in the general case one embeds $\g$
into some $\gl_n$ and uses Lemma~\ref{Lem:2.3}. It follows that only
finitely many subspaces of $\t$ are $G$-conjugate to a Cartan
subalgebra of a semisimple Lie algebra. Now it remains to note that
$G$-conjugate subspaces of $\t$ are $W$-conjugate. Here $W$ denotes
the Weyl group of $\g$.

{\it Step 2.} Conjugating $\widehat{\s}$ by an element of $Z_G(\s)$,
one may assume that there is a Cartan subalgebra $\t_0\subset
\widehat{\s}$ contained in $\t$. Since
$\widehat{\s}=[\widehat{\s},\widehat{\s}]+\s$, we see that $\t_0$ is
a sum of $\s$ and a Cartan subalgebra of a semisimple subalgebra of
$\g$. By step 1, there are only finitely many possibilities for
$\t_0$.

{\it Step 3.} Clearly, $\z(\widehat{\s})=\t_0\cap (\t_0\cap [\widehat{\s},\widehat{\s}])^\perp$,
where the orthogonal complement is taken with respect to some invariant non-degenerate
symmetric form on $\g$. Thus, by the previous steps, there are only finitely many possibilites
for $\z(\widehat{\s})$. Obviously, $\widehat{\s}$ is a direct sum of $\z(\widehat{\s})$
and a semisimple subalgebra of $\z_\g(\z(\widehat{\s}))$. Thence, $\widehat{\s}$ belongs to
one of finitely many $Z_G(\z(\widehat{\s}))$-conjugacy   classes of subalgebras.
To complete the proof of the lemma it remains to apply Lemma~\ref{Lem:2.3} to
$\g_0=\z_\g(\z(\widehat{\s}))$.
\end{proof}
\end{proof}

\begin{Cor}\label{Cor:2.5}
There are proper subspaces $\a_1,\ldots,\a_m\subset \a_{\g,\h}$ satisfying the following condition:
 if  $\widehat{\h}$ is a subalgebra of $\g$ containing $\h$ such that $c_{\g,\widehat{\h}}=c_{\g,\h}$
and $\a_{\g,\widehat{\h}}\not\subset \a_i$ for any $i$, then $\widehat{\h}\subset \h^{sat}$.
\end{Cor}
\begin{proof}
For $\a_i$ we take  elements of the set $\{\a_{\g,\widehat{\h}}|
\widehat{\h}=[\widehat{\h},\widehat{\h}]+R_u(\widehat{\h})+\h,\a_{\g,\widehat{\h}}\neq
\a_{\g,\h}\}$. Put
$\widehat{\h}_0=[\widehat{\h},\widehat{\h}]+R_u(\widehat{\h})+\h$.
Clearly,
$\widehat{\h}_0=[\widehat{\h}_0,\widehat{\h}_0]+R_u(\widehat{\h}_0)+\h$.
If $\a_{\g,\widehat{\h}_0}$ is not contained in any $\a_i$, then
$\a_{\g,\widehat{\h}_0}=\a_{\g,\h}$. Moreover, since $\h\subset
\widehat{\h}_0\subset \widehat{\h}$, we get
$c_{\g,\h}=c_{\g,\widehat{\h}}\leqslant
c_{\g,\widehat{\h}_0}\leqslant c_{\g,\h}$. Applying the following
lemma to $\g_0=\widehat{\h}_0,\h$, we get $\widehat{\h}_0=\h$.

\begin{Lem}\label{Lem:2.10}
For any algebraic subgroup $G_0\subset G$ we have
$$2(\dim \g-\dim\g_0)\geqslant 2c_{\g,\g_0}+2\dim \a_{\g,\g_0}+\dim\g-\dim\z_\g(\a_{\g,\g_0})$$
with the equality provided $G_0$ is observable.
\end{Lem}
\begin{proof}[Proof of Lemma \ref{Lem:2.10}]
 This follows from \cite{Knop1}, S\"atze 7.1,8.1,
Korollar 8.2\end{proof}

It follows that $\h$ is an ideal of $\widehat{\h}$ and that
$\widehat{\h}/\h$ is a commutative reductive algebraic Lie algebra.
Let $\widehat{H}$ denote the connected subgroup of $G$ corresponding
to $\widehat{\h}$. By Proposition 4.7 from \cite{Cartan},
$\widehat{H}/H^\circ$ acts on $G/H^\circ$ by central automorphisms,
equivalently, $\widehat{\h}\subset\h^{sat}$.
\end{proof}

The following lemma is used on step 3 of the proof of Proposition \ref{Prop:2.0}.

\begin{Lem}\label{Lem:2.7}
Let a highest weight $\lambda$ satisfy  condition (3) of Proposition
\ref{Prop:2.1}. Then \begin{itemize}\item[(3')] Any subgroup
$\widehat{H}\subset G$
  such that $H\subset \widehat{H}, H^\circ=\widehat{H}^\circ$ and
$V(\lambda)^H=V(\lambda)^{\widehat{H}}$ is contained in
$H^{sat}$.\end{itemize}
\end{Lem}
\begin{proof}
 By the Frobenius reciprocity,
$\C[G/\widehat{H}]^{(B)}_\lambda= \C[G/H]^{(B)}_\lambda$. By the
choice of $\lambda$,  $\C(G/H)^{B}=\C(G/\widehat{H})^{B}$.
Equivalently, $\C(G/B)^{H}=\C(G/B)^{\widehat{H}}$. Applying the main
theorem of the Galois theory to the field $\C(G/B)^{H^\circ}$, we
see that the images of $H/H^\circ, \widehat{H}/H^\circ$ in
$\Aut(\C(G/B)^{H^\circ})$ (or, equivalently,
$\Aut(\C(G/H^\circ)^{B})$) coincide.
 By assertion 1 of Lemma \ref{Lem:2.2}, $\widehat{H}/H^\circ=(H/H^\circ)\Gamma$, where
$\Gamma\subset \A_{G,H^\circ}$. Assertion 2 of Lemma \ref{Lem:2.2}
implies that $H$ is a normal subgroup in $\widehat{H}$. In virtue of
the natural inclusion $\C(G/H)^{B}\hookrightarrow \C(G/H^\circ)^B$,
the group $\widehat{H}/H$ acts trivially on $\C(G/H)^{B}$. It
remains to apply assertion 1 of Lemma \ref{Lem:2.2} one more time.
\end{proof}

Now we define subspaces $\a_1,\ldots,\a_k\subset \a_{\g,\h}$ and
sublattices $\X_1,\ldots,\X_l\subset\X_{G,H}$ satisfying the
assumptions of Proposition \ref{Prop:2.1}.

Suppose that $\A_{G,H}$ is a finite group. Take for  $\a_1,\ldots,
\a_k\subset \a_{\g,\h}$  subspaces found in Corollary \ref{Cor:2.5}.
Let $\A_{G,H}\cong \bigoplus_{i=1}^l \Z/p_i^{a_i}\Z$, where
$p_1,\ldots, p_l$ are distinct primes.  Take for $\X_i$ the lattice
$\X_{G,\widetilde{H}_i}$, where $\widetilde{H}_i$ denotes a unique
subgroup of $\widetilde{H}$ such that $\#\widetilde{H}_i/H=p_i$.
Clearly, $\widetilde{H}_i/H, i=\overline{1,l},$ are all minimal
proper subgroups of $\A_{G,H}$.

Now suppose that $\A_{G,H}$ is a one-dimensional torus. For
$\a_1,\ldots,\a_{k-1}$ we take subspaces  found in Corollary
\ref{Cor:2.5} and for  $\a_k$ take the subspace $\a_{\g,\h^{sat}}$.

Proposition \ref{Prop:2.0} follows from Proposition \ref{Prop:2.1},
Lemma \ref{Lem:2.7} and the following proposition.

\begin{Prop}\label{Prop:2.11}
Let $\lambda$ be a highest weight satisfying conditions (1),(2) of
Proposition \ref{Prop:2.1} for $\a_1,\ldots,\a_k,\X_1,\ldots,\X_l$
defined above and condition (3') of Lemma \ref{Lem:2.7} (or only
condition (1) if $c_{\g,\h}=0$). Then $\lambda$ has the properties
indicated in Proposition \ref{Prop:2.0}.
\end{Prop}
\begin{proof}
 Set $V:=V(\lambda)^H$.  By the
choice of $\lambda$, $\g_v=\h$ and $G_v\cap H^{sat}=H$ for $v\in V$
in general position.

First of all, we consider the case $c_{\g,\h}=0$. In this case
$H^{sat}=N_G(H)$ (this stems directly from Definition \ref{Def:1.1}
since $\dim \C(G/H)^{(B)}_{\lambda}=1$ for any $\lambda\in
\X_{G,H}$). Further, $N_G(H^\circ)/H^\circ$ is commutative and
thence $\widehat{H}\subset N_G(H)$ for any $\widehat{H}$ with
$\widehat{H}^\circ=H^\circ$. Thus $G_v\subset H^{sat}$ for a
non-zero vector $v\in V$. It follows from the choice of $\lambda$
that $G_v=H$.

In the sequel we assume that $c_{\g,\h}>0$.
Let us prove that the set
\begin{equation*}
\bigcup_{\widetilde{H}\supset H,
\widetilde{H}^\circ=H^\circ}V(\lambda)^{\widetilde{H}}\end{equation*}
is not dense in $V$. Any subgroup $\widetilde{H}\subset G$ with
$\widetilde{H}^\circ=H^\circ$ lies in $N_G(H^\circ)$. Denote by
$Y_n$ the subset of $N_G(H^\circ)/H^\circ$ consisting of all
elements $h$ such that $h$ and $H/H^\circ$ generate a finite
subgroup in $N_G(H)$, whose order divide $n$. For $h\in Y_n$ we
denote by $\widetilde{H}(h)$ the inverse image in $N_G(H^\circ)$ of
the subgroup of $N_G(H^\circ)/H^\circ$  generated by $h$ and
$H/H^\circ$.

Note that for every $n$ the subset $Y_n\subset N_G(H^\circ)/H^\circ$
is closed. Put $$Y_{n,i}=\{h\in Y_n| \codim_V
V(\lambda)^{\widetilde{H}(h)}=i\}.$$ This is a locally closed
subvariety of $Y_n$. Taking into account Lemma \ref{Lem:2.7}, we see
that $Y_{n,0}=\{1\}$ or $\varnothing$.

It is enough to show that for all $n,i>0$
the subset
\begin{equation}\label{eq:2.2}\bigcup_{h\in Y_{n,i}} V(\lambda)^{\widetilde{H}(h)}
\end{equation} is not dense in $V$.

Assume the converse: let $n,i\in \N$ be such that the subset (\ref{eq:2.2}) is dense
in $V$.
Then (compare with the proof of Proposition \ref{Prop:2.1})
$\dim Y_{n,i}\geqslant i$.
 It follows that $i\leqslant\dim Y_{n,i}\leqslant \dim G$. For $h_1,h_2\in Y_{n,i}$ the inequality
\begin{equation}\label{eq:2.3}\dim V(\lambda)^{\widetilde{H}(h_1)}\cap
V(\lambda)^{\widetilde{H}(h_2)}\geqslant \dim V-2i\geqslant \dim
V-2\dim G\end{equation} holds. Let $\widetilde{H}(h_1,h_2)$ denote
the algebraic subgroup of $G$ generated by
$\widetilde{H}(h_1),\widetilde{H}(h_2)$. Note that $\dim
V(\lambda)^{\widetilde{H}(h_1h_2)}=V(\lambda)^{\widetilde{H}(h_1)}\cap
V(\lambda)^{\widetilde{H}(h_2)}$. In virtue of (\ref{eq:2.3}) and
condition (2) of Proposition \ref{Prop:2.1},
$V(\lambda)^{\widetilde{H}(h_1,h_2)}\neq 0$.  By the choice of
$\lambda$, $\a_{\g,\widetilde{\h}(h_1,h_2)}=\a_{\g,\h}$. Therefore,
see Lemma \ref{Lem:2.10}, if $\widetilde{\h}(h_1,h_2)\neq \h$, then
$c_{\g,\widetilde{\h}(h_1,h_2)}<c_{\g,\h}$. But in this case
(\ref{eq:2.3}) contradicts condition (2) of Proposition
\ref{Prop:2.1}. So $\widetilde{\h}(h_1,h_2)=\h$ for any $h_1,h_2\in
Y_{n,i}$. In particular, any $h_1,h_2\in Y_{n,i}$  generate a finite
subgroup in $ N_G(H^\circ)/H^\circ$. Choose an irreducible component
$Y'\subset Y_{ni}$ of positive dimension.  Consider the map
$\rho:Y'\times Y'\rightarrow N_G(H^\circ)/H^\circ,(h_1,h_2)\mapsto
h_1 h_2^{-1}$. Its image is a non-discrete constructible set, whose
elements have  finite order in $N_G(H^\circ)/H^\circ$. Note that $1$
is a nonisolated point in $\overline{\operatorname{im}\rho}$. Thus
there is a locally closed subvariety $Z\subset
\overline{\operatorname{im}\rho}$ of positive dimension, whose
closure contains $1$. The subsets $Z_{j}:=\{z\in Z| z^j=1\}$ are
closed in $Z$. Thus $1\in \overline{Z}_j$ for some $j$. However, $1$
is an isolated point in $\{g\in N_G(H^\circ)/H^\circ| g^j=1\}$.
Contradiction.
\end{proof}

\section{Proof of Theorem \ref{Thm:1.3}}\label{SECTION_Thm2}
Again, one implication in Theorem \ref{Thm:1.3} is almost trivial.

\begin{proof}[Proof of $(a)\Rightarrow (b)$]
Let $V(\lambda)$ be such a simple module.  By Theorem \ref{Thm:1.2},
$\A_{G,H}$ is either a finite cyclic group or a one-dimensional
torus. Suppose that $\A_{G,H}\cong \C^\times$. As we noted in the
proof of the implication $(a)\Rightarrow (b)$, $\A_{G,H}$ acts on
$V(\lambda)^H$ by constants. If $\A_{G,H}\cong \C^\times$, then
$0\in\overline{\A_{G,H}v}$ for any $v\in V(\lambda)^H$. Thus
$0\in\overline{N_G(H)v}$ whence $0\in \overline{Gv}$. Contradiction.
\end{proof}

The proof of the other implication is much more complicated. Below
we assume that $\A_{G,H}$ is cyclic. At first, we prove
$(b)\Rightarrow (a)$ for reductive subgroups $H\subset G$ satisfying
the following condition.
\begin{itemize}
\item[(*)] The group $T_0:=(N_G(H)/H)^\circ$ is a torus, equivalently,
the  Lie algebra $\g^H$ is commutative.
\end{itemize}

The proof for $H$ satisfying (*) is based on the following technical
proposition, which is analogous to Proposition \ref{Prop:2.1}.

\begin{Prop}\label{Prop:3.0}
Let $H$ satisfy (*) and $\a_1,\ldots,\a_k,\X_1,\ldots,\X_l$ be such
as in Proposition \ref{Prop:2.1}. Then there is a highest weight
$\lambda$ satisfying conditions (1)-(3) of Proposition
\ref{Prop:2.1} (only (1) for $c_{\g,\h}=0$) and the following
condition:
\begin{itemize}
\item[(4)] The cone spanned by the weights of $T_0$ in $V(\lambda)^H$
coincides with the whole space $\X(T_0)\otimes_\Z\Q$.
\end{itemize}
\end{Prop}

We note that if $c_{\g,\h}=0$, then (4) holds automatically.

\begin{proof}[Proof of $(b)\Rightarrow (a)$ for $H$ satisfying (*)]
By Lemma \ref{Lem:2.7} and Proposition \ref{Prop:2.11}, there is a
dense subset $V^0\subset V:=V(\lambda)^H$ such that $G_v=H$ for any
$v\in V^0$. By condition(4) of Proposition \ref{Prop:3.0}, a general
orbit for the action $T_0:V$ is closed. It follows that there is
$v\in V$ such that $G_v=H$ and the orbit $N_G(H)v$ is closed. By the
Luna theorem, see \cite{VP}, Theorem 6.17, $Gv$ is also closed.
\end{proof}

In the proof of Proposition \ref{Prop:3.0} we will need several
lemmas. We may and will assume that $c_{\g,\h}>0$.

Let us introduce some notation. Set $L:=\X(T_0)\otimes_\Z\Q$. Let
$\Psi$ (resp., $\Psi^0$) denote the set of highest weights $\lambda$
with $V(\lambda)^H\neq 0$ (resp., satisfying condition (3)). By
Lemma \ref{Lem:2.1}, $\Psi$ is a monoid. For $\lambda\in\Psi$ by
$S(\lambda)$ we denote the set of weights of $T_0$ in
$V(\lambda)^H$. Since
$\C[G/H]^{(B)}_{\lambda^*}\C[G/H]^{(B)}_{\mu^*}\subset
\C[G/H]^{(B)}_{\lambda^*+\mu^*}$, we have $S(\lambda)+S(\mu)\subset
S(\lambda+\mu)$.  Finally, we denote by $\widetilde{H}$ the inverse
image of $T_0$ in $N_G(H)$ under the natural epimorphism
$N_G(H)\twoheadrightarrow N_G(H)/H$.

\begin{Lem}\label{Lem:3.1}
There is a highest weight $\nu$ satisfying conditions (1),(3),(4).
\end{Lem}
\begin{proof}
{\it Step 1.} Let us check that $\a_{\g,\widetilde{\h}}=\a_{\g,\h}$.
Since $\A_{G,H}$ is finite, Lemma \ref{Lem:2.2} implies that the
action $T_0:\C(G/H)^B$ is locally effective. It follows that
$c_{\g,\widehat{\h}}=c_{\g,\h}-\dim T_0$. The required equality
follows from the inclusion $\a_{\g,\widetilde{\h}}=\a_{\g,\h}$ and
Lemma \ref{Lem:2.10}.

{\it Step 2.} By step 1, elements $\lambda_0^*$ with $ \lambda_0\in
\Psi, 0\in S(\lambda_0),$ span $\a_{\g,\h}$. Clearly, $\Psi^0$ is an
ideal in $\Psi$. Therefore even $\lambda_0^*$ with
$\lambda_0\in\Psi_0:=\{\lambda_0\in \Psi^0| 0\in S(\lambda_0)\}$
span $\a_{\g,\h}$. Fix $\lambda_0\in \Psi_0$. We claim that
$S(\lambda_0)$ spans $L$. Indeed, otherwise there is a subgroup
$\widetilde{H}_0\subset \widetilde{H}$ such that
$\dim\widetilde{H}_0/H>0$ and $\widetilde{H}_0$ acts trivially on
$V(\lambda)^{H}$. By (3), $\widetilde{H}_0$ acts trivially on
$C(G/H)^{B}$, which contradicts $\#\A_{G,H}<\infty$.

{\it Step 3.} Set $\nu_0:=\lambda_0+\lambda_0^*$. Clearly,
$V(\lambda_0)^H\cong (V(\lambda_0^*)^H)^*$. Thus
$S(\lambda_0)=-S(\lambda_0^*)$.  It follows that $S(\nu_0)\supset
S(\lambda_0),-S(\lambda_0)$ whence the cone spanned by $S(\nu_0)$
coincides with $L$.

{\it Step 4.} Let $\mu',n$ be such as in Lemma \ref{Lem:3.0}. For
sufficiently large $m$ the cone spanned by $m S(\nu_0)+i S(\mu')$
coincides with $L$ for any $i=\overline{1,n}$. Thus for appropriate
$\mu'$ the weight $\nu:=m\nu_0+i\mu'$ satisfies (1),(3),(4).
\end{proof}

\begin{proof}[Proof of Proposition \ref{Prop:3.0}]
Let $\nu$ be such as in Lemma \ref{Lem:3.1}, $n$ be such as in Lemma
\ref{Lem:3.0}. We fix a norm $|\cdot|$ on $\a_{\g,\h}(\R)$ such that
$|\lambda|=|\lambda^*|$ for any $\lambda\in \a_{\g,\h}$. Let
$A_0,Y_0$ be such as in the proof of Proposition \ref{Prop:2.1}.

We choose $\lambda\in \Psi$ and  $A\in \R$  such that $\dim
V(\lambda)^H>A|\lambda|^{c-1}$, where $c:=c_{\g,\h}$, and $$
A|\lambda|^{c-1}>A_0(2|\lambda|+|\nu|n)^{c-1}+2\dim G+\dim Y_0.$$

By Lemma \ref{Lem:3.0} there is $i=\overline{1,n}$ such that
$\widetilde{\lambda}:=\lambda+\lambda^*+i\mu$ satisfies (1) and
automatically (3). As in the  proof of Proposition \ref{Prop:2.1},
$\widetilde{\lambda}$ satisfies (2). Finally, note that
$S(\lambda)=-S(\lambda^*)$. It follows that $S(\nu)\subset
S(\widetilde{\lambda})$ whence $\widetilde{\lambda}$ satisfies (4).
\end{proof}

\begin{proof}[Proof of Theorem \ref{Thm:1.3} in the general case]
Now $H$ is a subgroup of $G$ such that $\A_{G,H}$ is a finite cyclic
group and the algebra $\g^H$ is not commutative.

There is a finite cyclic subgroup $\Gamma\subset N_G(H)/H$ of  such
that
 $Z_{N_G(H)/H}(\Gamma)^\circ$ is a maximal torus of
$N_G(H)/H$, $\#\Gamma$ is prime and does not divide
$\#\A_{G,H},\#H/H^\circ$. Let $\overline{H}$ denote the inverse
image of $\Gamma$ in $N_G(H)$. Clearly, $\overline{H}\cap
H^{sat}=H$. Moreover, $(N_G(\overline{H})/\overline{H})^\circ$ is a
torus. Choose a highest weight  $\lambda$ satisfying conditions
(1)-(4) of Propositions \ref{Prop:2.1},\ref{Prop:3.0} (for
$\overline{H}$ instead of $H$). Let us check that $V(\lambda)$ has
the required properties.

Choose $\a_1,\ldots,\a_k,\X_1,\ldots,\X_l$ as in Proposition
\ref{Prop:2.11} for $\overline{H}$ instead of $H$. Let us check that
$\lambda$ satisfies conditions (1),(2) of Proposition \ref{Prop:2.1}
and condition (3') of Lemma \ref{Lem:2.7} for $H$. Condition (1) is
follows from the equality $\A_{G,H}=\A_{G,\overline{H}}$, which, in
turn, stems from \cite{Knop8}, Theorem 6.3, and the choice of
$\Gamma$. To check condition (2) it is enough to check that the
subset $Z\subset V(\lambda)$ defined there is closed. This will
follow if we check that $c_{\g,\widehat{\h}}<c_{\g,\h}$ for any
algebraic subalgebra $\widehat{\h}\subset\g$  such that
$\h\subsetneq \widehat{\h}, V(\lambda)^{\widehat{\h}}\neq\{0\}$. At
first, suppose that
$\widehat{\h}=[\widehat{\h},\widehat{\h}]+R_u(\widehat{\h})+\h$.
Then, by the choice of $\a_i$, we see that
$\a_{\g,\h}=\a_{\g,\widehat{\h}}$. Contradiction with Lemma
\ref{Lem:2.10}. Now let $\s$ denote a maximal reductive subalgebra
of $\widehat{\h}$ containing $\h$. Then $\s\subset
\s_0:=\h+\z(\s)\supsetneq \h$. It follows that
$c_{\g,\s_0}=c_{\g,\h}$. Thanks to Lemma \ref{Lem:2.2}, the last
equality contradicts $\#\A_{G,H}<\infty$. So conditions (1),(2) for
$\lambda$ and $H$ are checked.

Let us check condition (3'). Let $\widehat{H}$ be such a subgroup of
$G$ strictly containing $H$ such that $H^\circ=\widehat{H}^\circ,
V(\lambda)^H=V(\lambda)^{\widehat{H}}$. Let $\widetilde{H}$ denote
the algebraic subgroup of $G$ generated by
$\overline{H},\widehat{H}$. Then
$V(\lambda)^{\widetilde{H}}=V(\lambda)^{\overline{H}}\cap
V(\lambda)^{\widehat{H}}=V(\lambda)^{\overline{H}}$. Thanks to Lemma
\ref{Lem:2.7}, $\widetilde{H}\subset \overline{H}^{sat}$.  From the
choice of $\X_j$ it follows that $\widehat{H}\subset
\widetilde{H}=\overline{H}$. By the choice of $\Gamma$,
$\overline{H}=\widehat{H}$. So
$V(\lambda)^H=V(\lambda)^{\overline{H}}$. Choose a nilpotent element
$\xi\in\g^H$. Then $\exp(t\xi)\overline{H}\exp(t\xi)^{-1}\neq
\overline{H}$ but
$\exp(t\xi)V(\lambda)^{\overline{H}}=V(\lambda)^{\overline{H}}$.
But, by the proof of Proposition \ref{Prop:2.1}, there is $v\in
V(\lambda)^{\overline{H}}$ with $G_v=\overline{H}$. However,
$\exp(t\xi)v\not\in V(\lambda)^{\overline{H}}$. Contradiction. So
condition (3') holds for $\lambda,H$. By Proposition
\ref{Prop:2.11}, there is a dense open subset $V^0\subset
V(\lambda)^H$ such that $G_v=H$ for any $v\in V^0$.

It remains to prove that there is $v\in V^0$ with closed $G$-orbit
or, equivalently (by the Luna theorem, \cite{VP}, Theorem 6.17),
$N_G(H)$-orbit. Let $u\in V(\lambda)^{\overline{H}}$ be such that
$G_u=\overline{H}$ and $N_G(\overline{H})u$ is closed. Since
$\#\Gamma$ does not divide $\#H/H^\circ$, we have
$N_G(\overline{H})\subset N_G(H)$. By the Luna theorem, $N_G(H)u$ is
closed. Since there is a closed $N_G(H)$-orbit in $V(\lambda)^H$ of
dimension $\dim N_G(H)/H$, a general orbit is also closed.
\end{proof}

\section{Some examples}\label{SECTION_Examples}
In Introduction we have noted that  the group $\A_{G,H}$ can be
computed for any algebraic subgroup $H\subset G$. However, in
general, the computation algorithm is rather involved. In this
section we give examples  when the application of our theorems  is
easy.

\begin{Ex}
Let $H$ be a spherical observable subgroup of $G$ (the former means
that $G/H$ is spherical). In this case every automorphism of $G/H$
is central, so $\A_{G,H}=N_G(H)/H$. The classification of reductive
sherical subgroups is known and in this case groups $N_G(H)/H$ are
easy to compute. Note also that $G/H$ can be embedded to any module
$V(\lambda)$ provided $\lambda\not\in \X_{G,H^{sat}}$ for any
subgroup $\widetilde{H}\subset G$ containing $H$. For example, let
$G=\SL_{2n+1},H=\Sp_{2n}$. In this case $N_G(H)/H$ is a
one-dimensional torus. In fact, $G/H$ can be embedded into
$\bigwedge^3\C^{2n+1}$ provided $n\geqslant 3$.
\end{Ex}

\begin{Ex}
Let $H$ be a finite subgroup of $G$. It follows from results of
\cite{Knop8} that in this case $\A_{G,H}\cong Z(G)/Z(G)\cap H$. So
any homogeneous space $G/H$, where $Z(G)$ is a cyclic group or a
one-dimensional torus, can be embedded into a simple module as a
closed subvariety.
\end{Ex}


\bigskip

{\Small Department of Higher Algebra, Faculty of Mechanics and
Mathematics, Moscow State University.

E-mail address: ivanlosev@yandex.ru}
\end{document}